\documentclass{amsart}

\usepackage{amssymb}

\usepackage{amscd}

\newcommand\cD{{\mathcal{D}}}
\newcommand\cE{{\mathcal{E}}}
\newcommand\cF{{\mathcal{F}}}
\newcommand\cG{{\mathcal{G}}}
\newcommand\cH{{\mathcal{H}}}

\newcommand\cK{{\mathcal{K}}}
\newcommand\cL{{\mathcal{L}}}

\newcommand\cW{{\mathcal{W}}}

\newcommand\RR{{\mathbb{R}}}
\newcommand\CC{{\mathbb{C}}}
\newcommand\ZZ{{\mathbb{Z}}}

\newcommand\QQ{{\mathbb{Q}}}
\newcommand\TT{{\mathbb{T}}}

\newcommand{\End}{\operatorname{End}}
\newcommand{\supp}{\operatorname{supp}}

\newtheorem{theorem}{Theorem}

\theoremstyle{definition}

\newtheorem{example}[theorem]{Example}

\theoremstyle{remark}
\newtheorem{remark}[theorem]{Remark}

\begin{document}

\title[Classical and quantum dynamics in transverse geometry]{Classical and quantum dynamics in transverse geometry of Riemannian foliations}


\author{Yuri A. Kordyukov}
\address{Institute of Mathematics, Russian Academy of Sciences, Ufa}
\curraddr{} \email{yurikor@matem.anrb.ru}
\thanks{Supported by the Russian Foundation of Basic Research, grant 07-01-00081-a.}


\begin{abstract}
First, we survey some results on classical and quantum dynamical
systems associated with transverse Dirac operators on Riemannian
foliations. Then we illustrate these results by two examples of
Riemannian foliations: a foliation given by the fibers of a
fibration and a linear foliation on the two-dimensional torus.
\end{abstract}

\maketitle


\section*{Introduction}
In classical mechanics, a particle, moving on a compact manifold
$M$, is described by a point of the phase space, which is the
cotangent bundle $T^*M$ of $M$. The evolution of the particle in the
phase space is governed by the Hamilton equations of motion. In
particular, a Riemannian metric $g$ considered as a function on
$T^*M$ is the Hamiltonian of a free particle on $M$, and the
corresponding motion is given by the geodesic flow on $T^*M$.

In quantum mechanics, a particle on a compact manifold $M$ is
described by a function in the Hilbert space $L^2(M)$ called the
wave function. The evolution of the quantum particle is determined
by the Schr\"odinger equation. The Hamiltonian of a free quantum
particle, moving on a Riemannian manifold $(M,g)$, is the
Laplace-Beltrami operator $\Delta_g$ associated with the Riemannian
metric $g$.

It would be more convenient for us to consider Heisenberg's picture
of quantum mechanics, which deals with observable quantities and
their dynamics. For a free particle on a compact manifold $M$,
classical observables are real-valued functions on the phase space
$T^*M$, and quantum observables are pseudodifferential operators on
$M$ considered as (unbounded) self-ad\-jo\-int operators in
$L^2(M)$. The quantum evolution of pseudodifferential operators is
described by the Heisenberg equations of motion. The classical
evolution of their principal symbols is induced by the action of the
geodesic flow. The quantum and the classical evolutions are related
by the Egorov theorem.

The purpose of this paper is to discuss classical and quantum
dynamical systems for a particular class of Riemannian singular
spaces, namely, for the leaf space $M/\cF$ of a Riemannian foliation
$\cF$ on a compact manifold $M$. In this case, it is very natural to
use ideas and notions of noncommutative geometry by A. Connes
\cite{Co}. The conormal bundle $N^*\cF$ of $\cF$ carries a natural
foliation, $\cF_N$, so that the leaf space $N^*\cF/\cF_N$ of this
foliation can be naturally considered as the cotangent bundle of
$M/\cF$. Using constructions of noncommutative geometry, one can
associate to the singular space $N^*\cF/\cF_N$ some noncommutative
algebra, denoted in the paper by $C^{\infty}_{prop}(G_{{\mathcal
F}_N}, |T{\mathcal G}_N|^{1/2})$, which plays the role of an algebra
of smooth functions on this space, that is, of the algebra of
classical observables. The algebra of quantum observables is given
by the algebra of transverse pseudodifferential operators
$\Psi^{0,-\infty}(M,\cF)$ on $M$. We consider quantum Hamiltonians
defined by transverse Dirac operators. Then the associated classical
dynamics on $C^{\infty}_{prop}(G_{{\mathcal F}_N}, |T{\mathcal
G}_N|^{1/2})$ is induced by the action of the transverse geodesic
flow on $N^*\cF$, and the quantum evolution of transverse
pseudodifferential operators is related with the classical evolution
by a version of Egorov's theorem for transversally elliptic
operators stated in \cite{egorgeo,matrix-egorov}.

These results are closely related to the reduction theory for
quantum Hamiltonian systems with symmetry (see, for instance,
\cite{Gu-S82,Gu-Uribe,Gu-Uribe90,Tate99,Zelditch92} and references
therein).

We believe that the results discussed in the paper will play some
important role in further investigations of spectral theory of
transversally elliptic operators on foliated manifolds as well as in
the study of asymptotic spectral problems for elliptic operators on
foliated manifolds, such as adiabatic limits.

The paper is organized as follows. In Section~\ref{s:def}, we review
the results on classical and quantum dynamics in transverse geometry
of Riemannian foliations mentioned above. In Section~\ref{s:ex}, we
discuss two examples of Riemannian foliations: a foliation given by
the fibers of a fibration and a linear foliation on the
two-dimensional torus.

\section{Definitions and main results}\label{s:def}
Throughout in the paper, $(M,{\cF})$ is a compact foliated manifold,
$\operatorname{dim} M=n, \operatorname{dim} \cF=p, p+q=n$. We will
consider foliated charts $\varkappa: U\subset M\rightarrow I^p\times
I^q$ on $M$ with coordinates $(x,y)\in I^p\times I^q$ ($I$ is the
open interval $(0,1)$) such that the restriction of $\cF$ to $U$ is
given by the sets $y={\rm const}$. We will also use the following
notation: $T\cF$ is the tangent bundle of $\cF$; $Q=TM/T\cF$ is the
normal bundle of $\cF$; $N^*{\mathcal F}$ is the conormal bundle of
$\cF$.

\subsection{Classical phase space}
To define the cotangent bundle of the leaf space $M/\cF$ of the
foliation $(M,{\cF})$, one can proceed as follows. The conormal
bundle $N^*{\mathcal F}$ carries a natural foliation, ${\mathcal
F}_N$, called the linearized or the lifted foliation. One can define
it by constructing its foliated atlas. Given a foliated chart
$\varkappa: U\subset M\rightarrow I^p\times I^q$ on $M$ with
coordinates $(x,y)\in I^p\times I^q$, there is the corresponding
chart in $T^*M$ with coordinates written as $(x,y,\xi,\eta)\in
I^p\times I^q\times \RR^p\times \RR^q$. In these coordinates, the
restriction of the conormal bundle $N^*{\mathcal F}$ to $U$ is given
by the equation $\xi=0$. So we have a chart $\varkappa_n :
U_1\subset N^*{\mathcal F} \stackrel{\sim}{\longrightarrow}
I^p\times I^q\times \RR^q$ on $N^*{\mathcal  F}$ with the
coordinates $(x,y,\eta)\in I^p\times I^q\times \RR^q$. This chart is
a foliated chart on $N^*{\mathcal F}$ for the linearized foliation
${\mathcal F}_N$, and the restriction of ${\mathcal F}_N$ to $U_1$
is given by the level sets $y= {\rm const}, \eta={\rm const}$.

If the foliation $\cF$ is given by the fibers of a fibration $f:
M\to B$ over a compact manifold $B$, then the leaf space
$N^*{\mathcal F}/{\mathcal F}_N$ of the linearized foliation
${\mathcal F}_N$ coincides with the cotangent bundle $T^*B$ of $B$.
Therefore, for an arbitrary foliation $\cF$, the leaf space
$N^*{\mathcal F}/{\mathcal F}_N$ can be naturally considered as the
cotangent bundle of the leaf space $M/\cF$.

\begin{remark}\label{r:1}
This construction can be viewed a particular example of the
symplectic reduction in the sense of \cite[Chapter III, Section
14]{LM87} (see also \cite{Li75,Li77}).

Let $(X,\omega)$ be a symplectic manifold, and $Y$ a submanifold of
$X$ such that the $2$-form $\omega_Y$ induced by $\omega$ on $Y$ is
of constant rank. Let ${\mathcal F}_Y$ be the characteristic
foliation of $Y$ relative to $\omega_Y$ (that is, $T{\cF}_Y$ is the
skew-orthogonal complement of $TY$ in $TX$). If the foliation
${\mathcal F}_Y$ is simple, that is, it is given by the fibers of a
surjective submersion $p$ of $Y$ to a smooth manifold $Z$, then $Z$
has a unique symplectic form $\omega_Z$ such that
$p^*\omega_Z=\omega_Y$. The symplectic manifold $(Z, \omega_Z)$ is
said to be the reduced symplectic manifold associated with $Y$. In a
particular case when the submanifold $Y$ is the preimage of a point
under the momentum map associated with the Hamiltonian action of a
Lie group, the symplectic reduction associated with $Y$ is the
Mardsen-Weinstein symplectic reduction \cite{MW}.

In the case under consideration, one can consider the symplectic
reduction associated with the coisotropic submanifold
$Y=N^*{\mathcal F}$ in the symplectic manifold $X=T^*M$. The
corresponding characteristic foliation ${\mathcal F}_Y$ is the
linearized foliation ${\mathcal F}_N$, and the leaf space
$N^*{\cF}/{\mathcal F}_N$ plays a role of the reduced phase space
$Z$.
\end{remark}

In the general case, the leaf space $N^*{\cF}/{\mathcal F}_N$ is not
a smooth manifold. We will use the ideas of the noncommutative
geometry in the sense of A. Connes \cite{Co} and introduce a
noncommutative algebra, $C^{\infty}_{prop}(G_{{\mathcal
F}_N},|T{\mathcal G}_N|^{1/2})$, as a noncommutative analogue of an
algebra of smooth functions on $N^*{\cF}/{\mathcal F}_N$.

First, we recall several notions (for more details, see e.g.
\cite{survey} and references therein). Let $\gamma : [0.1]\to M$ be
a continuous leafwise path in $M$ with the initial point
$x=\gamma(0)$ and the final point $y=\gamma(1)$ and $T_0$ and $T_1$
arbitrary smooth submanifolds (possibly, with boundary), transversal
to the foliation, such that $x\in T_0$ and $y\in T_1$. Sliding along
the leaves of the foliation $\cF$ determines a diffeomorphism
$H_{T_0T_1}(\gamma)$ of a neighborhood of $x$ in $T_0$ to a
neighborhood of $y$ in $T_1$, called the holonomy map along
$\gamma$. The differential of $H_{T_0T_1}(\gamma)$ at $x$ gives rise
to a well-defined linear map $Q_x\to Q_y$, which is independent of
the choice of transversals $T_0$ and $T_1$. This map is called the
linear holonomy map and denoted by $dh_{\gamma}$. The adjoint of
$dh_{\gamma}$ yields a linear map $dh_{\gamma}^*:N^*\cF_y\to
N^*\cF_x$.

The holonomy groupoid $G$ of ${\cF}$ consists of
$\sim_h$-equivalence classes of continuous leafwise paths in $M$,
where we set $\gamma_1\sim_h \gamma_2$, if $\gamma_1$ and $\gamma_2$
have the same initial and final points and the same holonomy maps.
The set of units $G^{(0)}$ is the manifold $M$. The multiplication
in $G$ is given by the product of paths. The corresponding source
map $s:G\rightarrow M$ and range map $r:G\rightarrow M$ are given by
$s(\gamma)=\gamma(0)$ and $r(\gamma)=\gamma(1)$. Finally, the
diagonal map $\Delta:M\rightarrow G$ takes any $x\in M$ to the
element in $G$ given by the constant path $\gamma(t)=x, t\in [0,1]$.
To simplify the notation, we will identify $x\in M$ with
$\Delta(x)\in G$. For any $x\in M$ the map $s$ maps $G^x$ on the
leaf $L_x$ through $x$. The group $G^x_x$ coincides with the
holonomy group of $L_x$. The map $s:G^x\rightarrow L_x$ is the
covering map associated with the group $G^x_x$, called the holonomy
covering.

Let us introduce the groupoid $G_{{\mathcal F}_N}$ as the set of all
$(\gamma,\nu)\in G\times N^*{\mathcal F}$ such that
$r(\gamma)=\pi(\nu)$. The source map $s_N:G_{{\mathcal
F}_N}\rightarrow N^*{\mathcal F}$ and the range map
$r_N:G_{{\mathcal F}_N}\rightarrow N^*{\mathcal F}$ are defined as
$s_N(\gamma,\nu)=dh_{\gamma}^{*}(\nu)$ and $r_N(\gamma,\nu)=\nu$. We
have a map $\pi_G:G_{{\mathcal F}_N}\rightarrow G$ given by
$\pi_G(\gamma,\nu)=\gamma$. If $\cF$ is a Riemannian foliation,
$G_{{\mathcal F}_N}$ coincides with the holonomy groupoid of the
linearized foliation ${\mathcal F}_N$, but, for a general foliation,
these two groupoids may be different. Observe also that the leaf of
$\cF_N$ through $\nu\in N^*{\mathcal F}$ can be described as the set
of all points $dh_{\gamma}^{*}(\nu)\in N^*{\mathcal F}$, where
$\gamma\in G, r(\gamma)=\pi(\nu)$ (here $\pi :T^*M\to M$ is the
bundle map).

The groupoid $G_{{\mathcal F}_N}$ carries a natural codimension $q$
foliation $\cG_N$. The leaf of $\cG_N$ through a point
$(\gamma,\nu)\in G_{{\mathcal F}_N}$ is the set of all $(\gamma',
\nu')\in G_{{\mathcal F}_N}$ such that $\nu$ and $\nu'$ lie in the
same leaf in $\cF_N$. For any vector bundle $V$ on $M$, denote by
$|V|^{1/2}$ the associated half-density vector bundle. Let
$|T{\mathcal G}_N|^{1/2}$ be the line bundle of leafwise
half-densities on $G_{{\mathcal F}_N}$ with respect to the foliation
${\mathcal G}_N$. It is easy to see that
\[
|T{\mathcal G}_N|^{1/2}=r_N^*(|T{\mathcal F}_N|^{1/2})\otimes
s_N^*(|T{\mathcal F}_N|^{1/2}),
\]
where $s_N^*(|T{\mathcal F}_N|^{1/2})$ and $r_N^*(|T{\mathcal
F}_N|^{1/2})$ denote the lifts of the line bundle $|T{\mathcal
F}_N|^{1/2}$ of leafwise half-densities on $N^*{\mathcal F}$ via the
source and the range mappings $s_N$ and $r_N$ respectively.

Let $\varkappa_0: U\to I^p\times I^q, \varkappa_1: U'\to I^p\times
I^q$ be two foliated charts, let $T_0=\varkappa_0^{-1}(\{0\}\times
I^q)$, $T_1=\varkappa_1^{-1}(\{0\}\times I^q)$ be the corresponding
local transversals, and let $\pi_0 : T_0\to I^q$, $\pi_1 : T_1\to
I^q$ be the corresponding diffeomorphisms. The foliated charts
$\varkappa_0$ and $\varkappa_1$ are called compatible, if, for any
$m\in U$ and $m'\in U'$ with $\pi_0(m)=\pi_1(m')$, there is a
leafwise path $\gamma$ from $m$ to $m'$ such that the corresponding
holonomy diffeomorphism $H_{T_0T_1}(\gamma)$ considered as a local
diffeomorphism of $I^q$ is the identity in a neighborhood of
$\pi_0(m)=\pi_1(m')$.

For any pair of compatible foliated charts $\varkappa_0$ and
$\varkappa_1$ denote by $W(\varkappa_0,\varkappa_1)$ the subset in
$G$, consisting of all $\gamma\in G$ such that
$s(\gamma)=\varkappa_0^{-1}(x,y)\in U$ and
$r(\gamma)=\varkappa_1^{-1}(x',y)\in U'$ for some $(x,x',y)\in
I^p\times I^p\times I^q$ and the corresponding holonomy
diffeomorphism $H_{T_0T_1}(\gamma)$ considered as a local
diffeomorphism in $I^q$ is the identity in a neighborhood of $y\in
I^q$. There is a coordinate map
$\Gamma:W(\varkappa_0,\varkappa_1)\to I^p\times I^p\times I^q$,
which takes an element $\gamma\in W(\varkappa_0,\varkappa_1)$ as
above to the corresponding triple $(x,x',y)\in I^p\times I^p\times
I^q$. As shown in \cite{Co79}, the coordinate neighborhoods
$W(\varkappa_0,\varkappa_1)$ form an atlas of a $(2p+q)$-dimensional
manifold (in general, non-Hausdorff and non-paracompact) on $G$.
Moreover, the groupoid $G$ is a smooth groupoid.

Now let $\varkappa: U\subset M\rightarrow I^p\times I^q,
\varkappa_1: U' \subset M\rightarrow I^p\times I^q$, be two
compatible foliated charts on $M$. Then the corresponding foliated
charts $\varkappa_n: U_1\subset N^*{\mathcal  F}\rightarrow
I^p\times I^q \times \RR^q, (\varkappa_1)_n: U'_1 \subset
N^*{\mathcal F} \rightarrow I^p\times I^q \times \RR^q,$ are
compatible with respect to the foliation ${\mathcal  F}_N$. So they
define a foliated chart $V$ on the foliated manifold $(G_{{\mathcal
F}_N},{\mathcal G}_N)$ with the coordinates $(x,x',y,\eta) \in
I^p\times I^p\times I^q \times \RR^q$, and the restriction of
${\mathcal  G}_N$ to $V$ is given by the level sets $y= {\rm const},
\eta={\rm const}$.

A section $\sigma\in C^{\infty}(G_{{\mathcal F}_N}, |T{\mathcal
G}_N|^{1/2})$ is said to be properly supported, if the restriction
of the map $r_N:G_{\cF_N}\to N^*\cF$ to $\supp \sigma$ is a proper
map. One can introduce the structure of involutive algebra on the
space $C^{\infty}_{prop}(G_{{\mathcal F}_N}, |T{\mathcal
G}_N|^{1/2})$ of smooth, properly supported sections of $|T{\mathcal
G}_N|^{1/2}$ by the formulas
\begin{equation}\label{e:algebra}
\begin{aligned}
\sigma_1\ast
\sigma_2(\gamma,\nu)&=\int_{(\gamma_1,\nu_1)(\gamma_2,\nu_2)=(\gamma,\nu)}
\sigma_1(\gamma_1,\nu_1)\sigma_2(\gamma_2,\nu_2),\quad \gamma\in G_{{\mathcal F}_N},\\
\sigma^*(\gamma,\nu)&=\overline{\sigma((\gamma,\nu)^{-1})},\quad
\gamma\in G_{{\mathcal F}_N},
\end{aligned}
\end{equation}
where $\sigma, \sigma_1, \sigma_2\in C^{\infty}_{prop}(G_{{\mathcal
F}_N}, |T{\mathcal G}_N|^{1/2})$. The formula for $\sigma_1\ast
\sigma_2$ should be interpreted in the following way. The
composition $(\gamma,\nu)=(\gamma_1,\nu_1)(\gamma_2,\nu_2)$ is
defined if $\nu_2=dh_{\gamma_1}^{*}(\nu_1)$ and equals
$(\gamma,\nu)=(\gamma_1\gamma_2,\nu_1)$. Then we have
\begin{align*}
\sigma_1(\gamma_1,\nu_1)\sigma_2(\gamma_2,\nu_2) \in &
|T_{\nu_1}{\mathcal F}_N|^{1/2}\otimes
|T_{dh_{\gamma_1}^{*}(\nu_1)}{\mathcal F}_N|^{1/2}\otimes
|T_{\nu_2}{\mathcal F}_N|^{1/2}\otimes |T_{dh_{\gamma_2}^{*}(\nu_2)}{\mathcal F}_N|^{1/2} \\
 & \cong |T_{\nu_1}{\mathcal F_N}|^{1/2}\otimes |T_{\nu_2}{\mathcal F}_N|^{1}\otimes
 |T_{dh_{\gamma_1\gamma_2}^{*}(\nu_1)}{\mathcal F}_N|^{1/2},
\end{align*}
and, integrating the $|T_{\nu_2}{\mathcal F}_N|^{1}$-component of
$\sigma_1(\gamma_1,\nu_1)\sigma_2(\gamma_2,\nu_2)$ over the set of
all $\nu_2\in N^*\cF$ of the form $\nu_2=dh_{\gamma_1}^{*}(\nu)$,
which coincides with the leaf of $\cF_N$ through $\nu$, we get a
well-defined section of the bundle $r_N^*(|T{\mathcal
F}_N|^{1/2})\otimes s_N^*(|T{\mathcal F}_N|^{1/2})=|T\cG_N|^{1/2}. $

As mentioned above, the algebra $C^{\infty}_{prop}(G_{{\mathcal
F}_N}, |T{\mathcal G}_N|^{1/2})$ plays a role of the noncommutative
analogue of an algebra of functions on the leaf space
$N^*\cF/\cF_N$.

\subsection{Transversal pseudodifferential operators}\label{s:trpdo} Our quantum phase space
is a certain completion of the algebra $\Psi^{*,-\infty}(M,{\mathcal
F},E)$ of transversal pseudodifferential operators introduced in
\cite{noncom}. So we will briefly recall the definition of this
algebra. Let $E$ be a complex vector bundle over $M$ of rank $r$. We
will consider pseudodifferential operators, acting on
half-den\-si\-ti\-es. Let $C^{\infty}(M,E)$ denote the space of
smooth sections of the vector bundle $E\otimes |TM|^{1/2}$,
$L^2(M,E)$ the Hilbert space of square integrable sections of
$E\otimes |TM|^{1/2}$, ${\cD}'(M,E)$ the space of distributional
sections of $E\otimes |TM|^{1/2}$, ${\cD}'(M,E)=C^{\infty}(M,E)'$.
Finally, let $\Psi^{m}(M,E)$ denote the standard classes of
pseudodifferential operators, acting on $C^{\infty}(M,E)$.

For any $k_A \in S ^{m} (I ^{p} \times I^p\times I^q\times
{\RR}^{q}, {\cL}({\CC}^r))$, define an operator $A:
C^\infty_c(I^n,\CC^r)\to C^\infty(I^n,\CC^r)$ by the formula
\begin{equation}\label{loc}
Au(x,y)=(2\pi)^{-q} \int e^{i(y-y')\eta}k_A(x,x',y,\eta) u(x',y')
\,dx'\,dy'\,d\eta,
\end{equation}
where $u \in C^{\infty}_{c}(I^{n}, {\CC}^r), x \in I^{p}, y \in
I^{q}$. The function $k_A$ is called the complete symbol of $A$. As
usual, we will consider only classical (or polyhomogeneous) symbols,
that is, those symbols, which can be represented as an asymptotic
sum of homogeneous (in $\eta$) components.

Let $\varkappa: U\subset M\rightarrow I^p\times I^q, \varkappa_1: U'
\subset M\rightarrow I^p\times I^q$ be a pair of compatible foliated
charts on $M$ equipped with trivializations of the bundle $E$ over
them. Any operator $A$ of the form~(\ref{loc}) with the Schwartz
kernel, compactly supported in $I^n\times I^n$, determines an
operator $A:C^{\infty}_c(U,\left. E\right|_U)\to
C^{\infty}_c(U',\left. E\right|_{U'})$, which extends to an operator
in $C^\infty(M,E)$ in a trivial way. The resulting operator is
called an elementary operator of class $\Psi
^{m,-\infty}(M,{\mathcal F},E)$.

The  class $\Psi ^{m,-\infty}(M,{\mathcal F},E)$ consists of all
operators $A$ in $C^{\infty}(M,E)$, which can be represented in the
form
\[
A=\sum_{i} A_i + K,
\]
where $A_i$ are elementary operators of class $\Psi
^{m,-\infty}(M,{\mathcal F},E)$, corresponding to a pair
$\varkappa_i,\varkappa'_i$ of compatible foliated charts, $K\in \Psi
^{-\infty}(M,E)$.

The principal symbol of the operator $A$ given by (\ref{loc}) is the
leafwise half-density
\begin{equation}\label{k-principal}
\sigma(A)(x,x^\prime,y,\eta)=k_{A,m}(x,x^\prime,y,\eta)|dx|^{1/2}|dx^\prime|^{1/2},
\end{equation}
where $k_{A,m}$ is the degree $m$ homogeneous component of the
complete symbol $k_A$.

The global definition of the principal symbol is given as follows.
Let $\pi^*E$ denote the lift of the vector bundle $E$ to
$\tilde{T}^*M=T^*M\setminus 0$ via the bundle map $\pi:\tilde{T}^*M
\to M$. Put $\tilde{N}^*\cF=N^*\cF\setminus 0$. Denote by
$\tilde{\cF}_N$ the restriction of $\cF_N$ and by
$G_{\tilde{\cF}_N}$ the groupoid associated with $\tilde{\cF}_N$.
Let ${\mathcal L}(\pi^*E)$ be the vector bundle on
$G_{\tilde{\cF}_N}$, whose fiber at a point $(\gamma,\nu)\in
G_{\tilde{\cF}_N}$ is the space $\cL((\pi^*E)_{s_N(\gamma,\nu)},
(\pi^*E)_{r_N(\gamma,\nu)})$ of linear maps from
$(\pi^*E)_{s_N(\gamma,\nu)}$ to $(\pi^*E)_{r_N(\gamma,\nu)}$.
Consider the space $C^{\infty}_{prop}(G_{\tilde{\cF}_N}, {\mathcal
L}(\pi^*E)\otimes |T{\mathcal G}_N|^{1/2})$ of smooth, properly
supported sections of ${\mathcal L}(\pi^*E)\otimes |T{\mathcal
G}_N|^{1/2}$. One can introduce the structure of involutive algebra
on $C^{\infty}_{prop}(G_{\tilde{\cF}_N}, {\mathcal L}(\pi^*E)\otimes
|T{\mathcal G}_N|^{1/2})$ by formulas similar to \eqref{e:algebra}.
Let $S^{m}(G_{{\mathcal F}_N},{\mathcal L}(\pi^*E)\otimes
|T{\mathcal G}_N|^{1/2})$ be the space of all $s\in
C^{\infty}_{prop}(G_{\tilde{\cF}_N},{\mathcal L}(\pi^*E)\otimes
|T{\mathcal G}_N|^{1/2})$ homogeneous of degree $m$ with respect to
the action of $\RR$ given by the multiplication in the fibers of the
vector bundle $\pi_G:G_{{\mathcal F}_N}\rightarrow G $.

The principal symbol $\sigma(A)$ of an operator $A\in
\Psi^{m,-\infty}(M,{\mathcal F},E)$ given in local coordinates by
the formula (\ref{k-principal}) is globally defined as an element of
the space $S^{m}(G_{{\mathcal F}_N},{\mathcal L}(\pi^*E)\otimes
|T{\mathcal G}_N|^{1/2})$. Thus, we have the half-density principal
symbol mapping
\[
\sigma: \Psi^{m,-\infty}(M,{\mathcal F},E)\rightarrow
S^m(G_{{\mathcal F}_N},{\mathcal L}(\pi^*E)\otimes |T{\mathcal
G}_N|^{1/2}),
\]
which satisfies
\[
\sigma(AB)=\sigma(A)\sigma(B),\quad \sigma(A^*)=\sigma(A)^*
\]
for any $A\in \Psi^{m_1,-\infty}(M,{\mathcal F},E)$ and $B\in
\Psi^{m_2,-\infty}(M,{\mathcal F},E)$.

\subsection{Riemannian foliations}
We will consider a particular class of Hamiltonians associated with
a naive analogue of a Riemannian metric on the leaf space $M/\cF$.
It only exists for a particular class of foliations, called
Riemannian foliations.

Let $M$ be a compact manifold equipped with a foliation $\cF$, $g_M$
a Riemannian metric on $M$. Let $H=T{\cF}^{\bot}$ be the orthogonal
complement of the tangent bundle of ${\cF}$. So we have the direct
sum decomposition
\begin{equation}\label{e:decomp}
TM= T\cF\oplus H.
\end{equation}
There are natural isomorphisms $H\cong Q$ and $H^*\cong Q^*\cong
N^*\cF$. One can also decompose the metric $g_M$ into the sum of its
leafwise and transverse components:
\begin{equation}\label{e:g}
g_M=g_F+g_H.
\end{equation}

In a foliated chart with coordinates $(x,y)\in I^p\times I^q$ the
transverse part $g_H$ can be written as $
g_H=\sum_{\alpha\beta}g_{\alpha\beta}(x,y)\theta^\alpha\theta^\beta
$, where $\theta^\alpha\in H^*$ is the (unique) lift of $dy^\alpha$
under the projection $I^p\times I^q\to I^q$. The metric $g_M$ is
called bundle-like, if in any foliated chart the components
$g_{\alpha\beta}$ are independent of $x$,
$g_{\alpha\beta}(x,y)=g_{\alpha\beta}(y)$. Equivalently, one can say
that $g_M$ is bundle-like, if, for any leafwise continuous path
$\gamma$ from $x$ to $y$, the corresponding linear holonomy map
$dh_\gamma : Q_x\to Q_y$ is an isometry (see, for instance,
\cite{Molino,Re} for more details). The foliation $\cF$ is called
Riemannian, if there exists a bundle-like metric on $M$. In the
sequel, we will always assume that $\cF$ is a Riemannian foliation
and $g_M$ is a bundle-like metric.

Let $P_F$ (resp. $P_H$) denotes the orthogonal projection operator
of $TM$ onto $T\cF$ (resp. $H$). Denote by $\nabla^L$ the
Levi-Civita connection on $TM$ defined by $g_M$. The following
formulas define a connection $\nabla$ in $H$:
\begin{equation}\label{e:adapt}
\begin{aligned}
\nabla_XN &=P_H[X,N],\quad X\in C^\infty(M,T\cF),\quad N\in
C^\infty(M,H)\\ \nabla_XN&=P_H\nabla^L_XN,\quad X\in
C^\infty(M,H),\quad N\in C^\infty(M,H).
\end{aligned}
\end{equation}
Remark that the first identity in \eqref{e:adapt} yields a canonical
flat connection in $H$, defined along the leaves of $\cF$, which
exists for an arbitrary foliation and is called the Bott connection.
It turns out that $\nabla$ depends only on the transverse part of
the metric $g_M$ and preserves the inner product of $H$. It will be
called the transverse Levi-Civita connection.

Let $\omega_\cF$ denote the leafwise Riemannian volume form of
$\cF$. Let $f\in H_x$ and let $\tilde{f}\in C^\infty(M,H)$ denote
any infinitesimal transformation of $\cF$, which coincides with $f$
at $x$. By the definition of infinitesimal transformation, the local
flow generated by $\tilde{f}$ preserves the foliation and gives rise
to a well-defined action on $\Lambda^pT^*\cF$. The mean curvature
vector field $\tau\in C^\infty(M,H)$ of $\cF$ is defined by the
identity
\[
L_{\tilde{f}}\omega_\cF=g_M(\tau, \tilde{f}) \omega_\cF
\]
If $e_1,e_2,\ldots,e_p$ is a local orthonormal frame in $T\cF$, then
\[
\tau=\sum_{i=1}^pP_H(\nabla^L_{e_i}e_i).
\]

\subsection{Transverse Dirac operators}
For any $x\in M$, denote by $Cl(Q_x)$ the complex Clifford algebra
of $Q_x$. Recall that, relative to an orthonormal basis
$\{f_1,f_2,\ldots,f_q\}$ of $Q_x$, $Cl(Q_x)$ is the complex algebra
generated by $1$ and $f_1,f_2,\ldots,f_q$, satisfying the relations
\[
f_\alpha f_\beta+f_\beta f_\alpha=-2\delta_{\alpha\beta}, \quad
\alpha, \beta=1,2,\ldots,q.
\]

Consider the $\ZZ_2$-graded vector bundle $Cl(Q)$ over $M$ whose
fiber at $x\in M$ is $Cl(Q_x)$. This bundle is associated to the
principal $SO(q)$-bundle $O(Q)$ of oriented orthonormal frames in
$Q$, $Cl(Q)=O(Q)\times_{O(q)}Cl(\RR^q)$. Therefore, the transverse
Levi-Civita connection $\nabla$ induces a natural leafwise flat
connection $\nabla^{Cl(Q)}$ on $Cl(Q)$ which is compatible with the
multiplication and preserves the $\ZZ_2$-grading on $Cl(Q)$. If
$\{f_1, f_2, \ldots , f_q\}$ is a local orthonormal frame in $T^HM$,
and $\omega^\gamma_{\alpha\beta}$ is the coefficients of the
connection $\nabla$: $\nabla_{f_\alpha}f_\beta
=\sum_\gamma\omega^\gamma_{\alpha\beta}f_\gamma$, then
\begin{equation}\label{e:spin}
\nabla^{Cl(Q)}_{f_\alpha} = f_\alpha+\frac{1}{4}\sum_{\gamma=1}^q
\omega^\gamma_{\alpha\beta}c(f_\beta)c(f_\gamma),
\end{equation}
where $c(a)$ denotes the action of an element $a\in
C^\infty(M,Cl(Q))$ on $C^\infty(M,Cl(Q))$ by pointwise left
multiplication.

A transverse Clifford module is a complex vector bundle $\cE$ on $M$
endowed with a fiberwise action of the bundle $Cl(Q)$. We will
denote the action of $a\in Cl(Q_x)$ on $s\in \cE_x$ as $c(a)s\in
\cE_x$. A transverse Clifford module $\cE$ is called self-adjoint if
it endowed with a Hermitian metric such that the operator $c(f) :
\cE_x\to \cE_x$ is skew-adjoint for any $x\in M$ and $f\in Q_x$. Any
transverse Clifford module $\cE$ carries a natural $\ZZ_2$-grading
$\cE=\cE_+\oplus\cE_-$.

A Hermitian connection $\nabla^\cE$ on a transverse Clifford module
$\cE$ is called a Clifford connection if it is compatible with the
Clifford action, that is, for any $f\in C^\infty(M,H)$ and $a\in
C^\infty(M,Cl(Q))$,
\[
[\nabla^{\cE}_f, c(a)]=c(\nabla^{Cl(Q)}_fa).
\]

A self-adjoint transverse Clifford module $\cE$ equipped with a
leafwise flat Clifford connection $\nabla^\cE$ is called a
transverse Clifford bundle. For any transverse Clifford bundle
$(\cE, \nabla^\cE)$, we define the operator $D^\prime_\cE$ acting on
the sections of $\cE$ as the composition
\[
C^\infty(M,\cE)\stackrel{\nabla^{\cE}}{\longrightarrow}
C^\infty(M,Q^*\otimes \cE) = C^\infty(M,Q\otimes \cE)
\stackrel{c}{\longrightarrow} C^\infty(M, \cE).
\]
Here we identify the bundle $Q$ and $Q^*$ by means of the metric
$g_M$. The transverse Dirac operator $D_\cE$ is defined as
\[
D_\cE=D^\prime_\cE-\frac12 c(\tau).
\]
If $f_1,\ldots,f_q$ is a local orthonormal frame for $H$, then
\[
D_\cE=\sum_{\alpha=1}^q
c(f_\alpha)\left(\nabla^{\cE}_{f_\alpha}-\frac12 g_M(\tau, f_\alpha)
\right).
\]
This operator is odd with respect to the natural $\ZZ_2$-grading on
$\cE$.

Denote by $(\cdot,\cdot)_x$ the inner product in the fiber $\cE_x$
over $x\in M$. Then the inner product in $L^2(M, \cE)$ is given by
the formula
\[
(s_1, s_2)=\int_M (s_1(x),s_2(x))_x\omega_M, \quad s_1, s_2\in
L^2(M, \cE),
\]
where $\omega_M=\sqrt{\det g}\,dx$ denotes the Riemannian volume
form on $M$. As shown in \cite{matrix-egorov}, the transverse Dirac
operator $D_\cE$ is formally self-adjoint in $L^2(M, \cE)$.

The transverse Dirac operators were introduced in
\cite{G-K91a,G-K91} (see \cite{matrix-egorov,vanishing} for further
references).

We will use the Riemannian volume form $\omega_M$ to identify the
half-densities bundle with the trivial one. So the action of $D_\cE$
on half-densities is defined by
\[
D_\cE (u|\omega_M|^{1/2})=(D_\cE u)|\omega_M|^{1/2}, \quad u\in
C^\infty(M,\cE).
\]

\begin{example}
Assume that $\cF$ is transversely oriented and the normal bundle $Q$
is spin. Thus the $SO(q)$ bundle $O(Q)$ of oriented orthonormal
frames in $Q$ can be lifted to a $Spin(q)$ bundle $O'(Q)$ so that
the projection $O'(Q)\to O(Q)$ induces the covering projection
$Spin(q)\to SO(q)$ on each fiber.

Let $F(Q), F_+(Q), F_-(Q)$ be the bundles of spinors
\[
F(Q)=O'(Q)\times_{Spin(q)}S, \quad
F_\pm(Q)=O'(Q)\times_{Spin(q)}S_\pm .
\]
Since $\dim Q=q$ is even $\End F(Q)$ is as a bundle of algebras over
$M$ isomorphic to the Clifford bundle $Cl(Q)$. The transverse
Levi-Civita connection $\nabla$ lifts to a leafwise flat Clifford
connection $\nabla^{F(Q)}$ on $F(Q)$. So $F(Q)$ is a transverse
Clifford bundle.

More generally, one can take a Hermitian vector bundle $\cW$
equipped with a leafwise flat Hermitian connection $\nabla^\cW$.
Then $F(Q)\otimes \cW$ is a transverse Clifford bundle: the action
of $a\in C^\infty(M,Cl(Q))$ on $C^\infty(M, F(Q)\otimes \cW)$ is
given by $c(a)\otimes 1$ ($c(a)$ denotes the action of $a$ on
$C^\infty(M, F(Q))$), and the product connection
$\nabla^{F(Q)\otimes \cW}=\nabla^{F(Q)}\otimes 1 + 1\otimes
\nabla^{\cW}$ on $F(Q)\otimes \cW$ is a Clifford connection.
\end{example}

\begin{example}
Another example of a transverse Clifford bundle associated with a
transverse almost complex structure on $(M,\cF)$, a transverse
Clifford module $\Lambda^{0,*}$, is described in \cite{vanishing}.
\end{example}

\begin{example}\label{ex:sign}
The decomposition (\ref{e:decomp}) induces a bigrading on $\Lambda
T^{*}M$:
\[
\Lambda T^{*}M=\bigoplus_{i,j=1}^{n}\Lambda^{i,j}T^{*}M, \quad
\Lambda^{i,j}T^{*}M=\Lambda^{i}T\cF^{*}\otimes \Lambda^{j} H^{*}.
\]
In this bigrading, the de Rham differential $d$ can be written as
\begin{equation}\label{e:d}
d=d_F+d_H+\theta,
\end{equation}
where $d_F$ and $d_H$ are first order differential operators (called
the tangential de Rham differential and the transversal de Rham
differential accordingly), and $\theta$ is a zero order differential
operator.

By definition, the transverse signature operator is a first order
differential operator in  $C^{\infty}(M,\Lambda H^{*})$ given by
\[
D_H=d_H + d^*_H.
\]
Consider the vector bundle $\cE=\Lambda H^*\otimes \CC$ equipped
with a natural structure of a transverse Clifford bundle. As shown
in \cite{matrix-egorov}, for the associated transverse Dirac
operator $D_{\cE}$, we have
\[
D_{\cE}= d_H+d^*_H-\frac12(\varepsilon_{\tau^*}+i_{\tau}),
\]
where $\varepsilon_{\tau^*}$ denotes the exterior product by the
covector $\tau^*\in Q_x^*$ dual to $\tau$, $i_\tau$ the interior
product by $\tau$. So we see that the transverse signature operator
$D_H$ coincides with $D_{\cE}$ if and only if $\tau=0$, that is, all
the leaves are minimal submanifolds.
\end{example}

\subsection{Quantum dynamics}
Let $D_\cE$ be a transverse Dirac operator associated with a
transverse Clifford bundle $(\cE, \nabla^\cE)$. The operator $D_\cE$
is essentially self-adjoint with initial domain $C^\infty(M,\cE)$.
Define an unbounded linear operator $\langle D_\cE\rangle$ in the
space $L^2(M,\cE)$ as
\[
\langle D_\cE\rangle=(D_\cE^2+I)^{1/2}.
\]
By the spectral theorem, the operator $\langle D_\cE\rangle $ is
well-defined as a positive, self-adjoint operator in $L^2(M,\cE)$.
It can be shown that the Sobolev space $H^1(M,\cE)$ is contained in
the domain of $\langle D_\cE\rangle$ in $L^2(M,\cE)$.

By the spectral theorem, the operator $\langle D_\cE\rangle$ defines
a strongly continuous group $e^{it\langle D_\cE\rangle}$ of bounded
operators in $L^2(M,\cE)$. Consider a one-parameter group $\Phi_t$
of $\ast$-automorphisms of the algebra ${\mathcal L}(L^2(M,\cE))$
defined by
\[
\Phi_t(T)=e^{i t\langle D_\cE\rangle}Te^{-i t\langle D_\cE\rangle},
\quad T\in {\mathcal L}(L^2(M,\cE)), \quad t\in \RR.
\]

By the spectral theorem, the operator $\langle
D_\cE\rangle^s=(D_\cE^2+I)^{s/2}$ is a well-defined positive
self-adjoint operator in $\cH=L^2(M,\cE)$ for any $s\in\RR$, which
is unbounded if $s>0$. For any $s\geq 0$, denote by $\cH^s$ the
domain of $\langle D_\cE\rangle^s$, and, for $s<0$,
$\cH^s=(\cH^{-s})^*$. Put also $\cH^{\infty}=\bigcap_{s\geq 0}\cH^s,
\quad \cH^{-\infty}=(\cH^{\infty})^*$. It is clear that $H^s(M,\cE)
\subset \cH^s$ for any $s\geqslant 0$ and $\cH^s \subset H^s(M,\cE)$
for any $s<0$.  In particular, $C^\infty(M,\cE) \subset \cH^s$ for
any $s$.

We say that a bounded operator $A$ acting on $\cH^{\infty}$ belongs
to $\cL(\cH^{-\infty},\cH^{\infty})$ (resp.
$\cK(\cH^{-\infty},\cH^{\infty})$), if, for any $s$ and $r$, it
extends to a bounded (resp. compact) operator from $\cH^s$ to
$\cH^r$, or, equivalently, the operator $\langle
D_\cE\rangle^rA\langle D_\cE\rangle^{-s}$ extends to a bounded
(resp. compact) operator on $L^2(M,\cE)$. It is easy to see that
$\cL(\cH^{-\infty},\cH^{\infty})$ is an involutive subalgebra in
$\cL(\cH)$ and $\cK(\cH^{-\infty},\cH^{\infty})$ is its ideal. We
also introduce the class $\cL^1(\cH^{-\infty},\cH^{\infty})$, which
consists of all operators from $\cK(\cH^{-\infty},\cH^{\infty})$
such that, for any $s$ and $r$, the operator $\langle
D_\cE\rangle^rA\langle D_\cE\rangle^{-s}$ is a trace class operator
on $L^2(M,\cE)$. It should be noted that any operator $K$ with
smooth kernel belongs to $\cL^1(\cH^{-\infty},\cH^{\infty})$.

\begin{theorem}\cite{matrix-egorov}
\label{Egorov1} Let $D_\cE$ be a transverse Dirac operator
associated with a transverse Clifford bundle $(\cE, \nabla^\cE)$.
For any $K\in \Psi^{m,-\infty}(M,{\mathcal F},\cE)$, there exists an
operator $K(t)\in\Psi^{m,-\infty}(M,{\mathcal F},\cE)$ such that
$\Phi_t(K)-K(t), t\in \RR,$ is a smooth family of operators of class
$\cL^1(\cH^{-\infty},\cH^{\infty})$.
\end{theorem}

By this theorem, we have a well-defined one-parameter automorphism
group
\[
\Phi_t : \hat\Psi^{0,-\infty}(M,{\mathcal F},\cE)\to
\hat\Psi^{0,-\infty}(M,{\mathcal F},\cE),\quad t\in \RR,
\]
of the completed algebra of transverse pseudodifferential operators
defined by
\[
\hat\Psi^{0,-\infty}(M,{\mathcal
F},\cE)=\Psi^{0,-\infty}(M,{\mathcal F},\cE)+
\cL^1(\cH^{-\infty},\cH^{\infty}).
\]
It describes the quantum dynamics associated with the quantum
Hamiltonian $\langle D_\cE\rangle$.

\subsection{Classical dynamics}
In this Section, we give a definition of the transverse geodesic
flow associated with a bundle-like metric $g$ on a compact foliated
manifold $(M,\cF)$. This notion is a particular case of the
bicharacteristic flow associated with a first order transversally
elliptic operator (see \cite{egorgeo,matrix-egorov}).

Denote by $f_t$ the geodesic flow of the Riemannian metric $g$ on
$T^*M$. Since $g$ is bundle-like, the flow $f_t$ preserves the
conormal bundle $N^*{\mathcal F}$, and its restriction to
$N^*{\mathcal F}$ (denoted also by $f_t$) preserves the foliation
$\cF_N$, that is, takes any leaf of $\cF_N$ to a leaf. Moreover, one
can show existence of a flow $F_t$ on $G_{{\mathcal F}_N}$ such that
$s_N\circ F_t=f_t\circ s_N$, $r_N\circ F_t=f_t\circ r_N$, which
preserves the foliation $\cG_N$.

Let $X$ be the generator of the geodesic flow $f_t$. Then it is
tangent to $N^*\cF$ and determines a vector field on $N^*\cF$,
denoted also by $X$. Since $f_t$ preserves the foliation $\cF_N$,
$X$ is an infinitesimal transformation of $\cF_N$, and there exists
a vector field $\cH$ on $G_{{\mathcal F}_N}$ such that $ds_N(\cH)=X$
and $dr_N(\cH)=X$. Then the vector field $\cH$ is the generator of
the flow $F_t$.

In a foliation chart, $X$ is given by
\[
X = \sum_{j=1}^p \left(\partial_{\xi_j} \sqrt{g} \frac{\partial
}{\partial x_j} - \partial_{x_j} \sqrt{g} \frac{\partial }{\partial
\xi_j}\right) + \sum_{k=1}^q \left(\partial_{\eta_k}
\sqrt{g}\frac{\partial }{\partial y_k} - \partial_{y_k} \sqrt{g}
\frac{\partial }{\partial \eta_k}\right),
\]
and $\cH$ is given by
\begin{multline*}
\cH(x,x',y,\eta) = \sum_{j=1}^p \partial_{\xi_j}
\sqrt{g}(x,y,0,\eta) \frac{\partial }{\partial x_j} +\sum_{j=1}^p
\partial_{\xi_j} \sqrt{g}(x^\prime,y,0,\eta) \frac{\partial
}{\partial x^\prime_j}
\\ + \sum_{k=1}^q \left(\partial_{\eta_k} \sqrt{g}(y,\eta)\frac{\partial
}{\partial y_k} -
\partial_{y_k} \sqrt{g}(y,\eta) \frac{\partial }{\partial
\eta_k}\right),\\ \quad (x,x',y,\eta) \in I^p\times I^p\times
I^q\times \RR^q.
\end{multline*}

The flow $F_t^*$ on $C^{\infty}_{prop}(G_{{\mathcal F}_N},
|T{\mathcal G}_N|^{1/2})$ induced by $F_t$ is called the transverse
geodesic flow.

Now consider a transverse Clifford bundle $(\cE, \nabla^\cE)$. The
pull-back of the connection form of $\nabla^\cE$ by the projection
$\pi: N^*\cF \to M$ yields a connection $\nabla^{\pi^*\cE}$ on
$\pi^*\cE$. The parallel transport along the orbits of the
transverse geodesic flow $f_t$ on $N^*\cF$ associated with the
connection $\nabla^{\pi^*\cE}$ defines a one-parameter automorphism
group $\alpha_t$ of the vector bundle $\pi^*\cE$, covering the flow
$f_t$, and the induced flow $\alpha_t^*$ on
$C^\infty(N^*\cF,\pi^*\cE)$. In its turn, this flow induces the flow
$\operatorname{Ad}(\alpha_t)^*$ on the space
$C^{\infty}_{prop}(G_{{\mathcal F}_N}, {\mathcal L}(\pi^*\cE)\otimes
|T{\mathcal G}_N|^{1/2})$, which satisfies
\[
\frac{d}{dt}\operatorname{Ad}(\alpha_t)^*\sigma=\nabla_{\cH}\sigma,
\quad \sigma\in C^{\infty}_{prop}(G_{{\mathcal F}_N}, {\mathcal
L}(\pi^*\cE)\otimes |T{\mathcal G}_N|^{1/2}).
\]
This flow will be called the transverse parallel transport. One can
show that
\[
\operatorname{Ad}(\alpha_t)^*\circ s^*_N=s^*_N\circ \alpha_t^*,
\quad \operatorname{Ad}(\alpha_t)^*\circ r^*_N=r^*_N\circ
\alpha_t^*.
\]
The following theorem proved in \cite{matrix-egorov} (for scalar
operators, see \cite{egorgeo}) relates the quantum dynamics defined
by a transverse Dirac operator with the corresponding classical
dynamics.

\begin{theorem}
\label{Egorov2} Let $D_\cE$ be a transverse Dirac operator
associated with a transverse Clifford bundle $(\cE, \nabla^\cE)$.
For any $K\in \Psi^{m,-\infty}(M,{\mathcal F},\cE)$ with the
principal symbol $\sigma\in S^m(G_{{\mathcal F}_N},{\mathcal
L}(\pi^*\cE)\otimes |T{\mathcal G}_N|^{1/2})$, the principal symbol
$\sigma_t\in S^m(G_{{\mathcal F}_N},{\mathcal L}(\pi^*\cE)\otimes
|T{\mathcal G}_N|^{1/2})$ of the operator $K(t)$ is given by
\[
\sigma_t=\operatorname{Ad}(\alpha_t)^*\sigma.
\]
\end{theorem}

\begin{remark}
The construction of the transversal geodesic flow provides an
example of what can be called noncommutative symplectic (or, maybe,
better, Poisson) reduction in the setting discussed above.

As in Remark~\ref{r:1}, let $(X,\omega)$ be a symplectic manifold,
and $Y$ a submanifold of $X$ such that the $2$-form $\omega_Y$
induced by $\omega$ on $Y$ is of constant rank. Let ${\mathcal F}_Y$
be the characteristic foliation of $Y$ relative to $\omega_Y$.
Suppose that the foliation ${\mathcal F}_Y$ is simple and $(B,
\omega_B)$ is the associated reduced symplectic manifold. If $Y$ is
invariant under the Hamiltonian flow of a Hamiltonian $H\in
C^\infty(X)$ (this is equivalent to the fact that
$(dH)\left|_Y\right.$ is constant along the leaves of the
characteristic foliation ${\mathcal F}_Y$), there exists a unique
function $\hat{H}\in C^\infty(B)$, called the reduced Hamiltonian,
such that $H\left|_Y\right. = \hat{H}\circ p$. Furthermore, the map
$p$ projects the restriction of the Hamiltonian flow of $H$ to $Y$
to the reduced Hamiltonian flow on $B$ defined by the reduced
Hamiltonian $\hat{H}$ (see, for instance, \cite[Chapter III, Theorem
14.6]{LM87}).

In the case under consideration, one can consider as above the
symplectic reduction associated with the coisotropic submanifold
$Y=N^*{\mathcal F}$. The Hamiltonian of the geodesic flow on $T^*M$
determined by a bundle-like metric $g_M$ on $M$ satisfies the above
invariance condition, and, if the foliation $\cF$ is simple, one can
apply the symplectic reduction procedure mentioned above, that gives
us the geodesic flow on the cotangent bundle of the base as the
corresponding reduced Hamiltonian flow (see also
Section~\ref{s:Riem-sub}). For a general foliation, following the
ideas of noncommutative geometry, one can consider the corresponding
reduced Hamiltonian flow as the flow $F_t^*$ on the algebra
$C^{\infty}_{prop}(G_{{\mathcal F}_N},|T{\mathcal G}_N|^{1/2})$.
Following the ideas of \cite{Block-Ge, Xu}, one can interpret the
algebra $C^{\infty}_{prop}(G_{{\mathcal F}_N},|T{\mathcal
G}_N|^{1/2})$ as a noncommutative Poisson manifold and the flow
$F_t^*$ as a noncommutative Hamiltonian flow.
\end{remark}

\section{Examples}\label{s:ex}
\subsection{Riemannian submersions}\label{s:Riem-sub} Suppose that a foliation $\cF$ on
a compact manifold $M$ is given by the fibers of a fibration $p:
M\to B$ over a compact manifold $B$. Then, for any $x\in M$,
$N^*_x\cF$ coincides with the image of the cotangent map $dp(x)^* :
T_{p(x)}^*B \to T_x^*M$. The inverse maps $(dp(x)^*)^{-1} : N_x^*\cF
\to T_{p(x)}^*B$ determine a fibration $N^*\cF \to T^*B$ whose
fibers are the leaves of the linearized foliation $\cF_N$. Thus,
$N^*\cF$ is diffeomorphic to the fiber product
\[
M\times_B T^*B=\{ (x,\xi)\in M\times T^*B : p(x)=\pi(\xi)\}
\]
with a diffeomorphism $M\times_B T^*B
\stackrel{\cong}{\longrightarrow} N^*\cF$, given by
\begin{equation}\label{e:nf}
(x,\xi)\in M\times_B T^*B\mapsto dp(x)^*(\xi)\in N^*_x\cF.
\end{equation}

The holonomy groupoid $G$ of $\cF$ is the fiber product
\[
M\times_BM=\{(x,y)\in M\times M : p(x)=p(y)\},
\]
where $s(x,y)=y, r(x,y)=x$. Similarly, the holonomy groupoid
$G_{\cF_N}$ is the fiber product $N^*\cF\times_{T^*B}N^*\cF$, which
consists of all $(\nu_x,\nu_y)\in N^*_x\cF \times N^*_y\cF$ such
that $(x,y)\in M\times_BM$ and $(dp(x)^*)^{-1}(\nu_x)
=(dp(x)^*)^{-1}(\nu_y)$, with $s_N(\nu_x,\nu_y)=\nu_y,
r_N(\nu_x,\nu_y)=\nu_x$. On the other hand,
$N^*\cF\times_{T^*B}N^*\cF$ is also diffeomorphic to the fiber
product
\[
G\times_B T^*B = \{(x,y,\xi)\in M \times M\times T^*B :
p(x)=p(y)=\pi_B(\xi)\}.
\]
A diffeomorphism $G\times_B T^*B \stackrel{\cong}{\longrightarrow}
N^*\cF\times_{T^*B}N^*\cF$ can be defined as
\[
(x,y,\xi)\in G\times_B T^*B \mapsto (dp(x)^*(\xi),dp(y)^*(\xi))\in
N^*\cF\times_{T^*B}N^*\cF.
\]
The foliation $\cG_N$ is given by the fibers of the fibration
\[
(x,y,\xi)\in G_{\cF_N} \cong G\times_B T^*B \mapsto \xi\in T^*B.
\]

As in Remark~\ref{r:1}, consider the cotangent bundle $T^*M$ as a
symplectic manifold equipped with the canonical symplectic
structure, and $N^*\cF$ as its closed coisotropic submanifold. Then
the linearized foliation $\cF_N$ coincides with the characteristic
foliation of this coisotropic submanifold. It is well-known (see,
for instance, \cite{GS79}) that the fiber product
$N^*\cF\times_{T^*B}N^*\cF$ is a canonical relation in ${T}^*M$,
which is often called the flowout of the coisotropic submanifold
${N}^*{\cF}$. For a complex vector bundle $E$ on $M$, the algebra of
Fourier integral operators in $C^\infty(M,E)$ associated with this
canonical relation coincides with the algebra
$\Psi^{*,-\infty}(M,{\cF},E)$. It is a particular case of a general
construction described in \cite{GS79}.

In this case, the corresponding classes of symbols
$S^{m}(G_{{\mathcal F}_N}, {\mathcal L}(\pi^*E)\otimes |T{\mathcal
G}_N|^{1/2})$ can be described as follows. For any $\xi\in T^*B$,
let $\Psi^{-\infty}((N^*\cF)_\xi, (\pi^*E)_\xi)$ be the involutive
algebra of all smoothing operators, acting on
$C^\infty((N^*\cF)_\xi, (\pi^*E)_\xi)$, where $(N^*\cF)_\xi$ is the
fiber of the fibration $N^*\cF\to T^*B$ at $\xi$ and $(\pi^*E)_\xi$
is the restriction of $\pi^*E$ to $(N^*\cF)_\xi$. Consider a field
$\Psi^{-\infty}(N^*\cF,\pi^*E)$ of involutive algebras on $T^*B$
whose fiber at $\xi\in T^*B$ is $\Psi^{-\infty}((N^*\cF)_\xi,
(\pi^*E)_\xi)$. For any section $K$ of the field
$\Psi^{-\infty}(N^*\cF,\pi^*E)$, the Schwartz kernels of the
operators $K_\xi$ in $C^\infty((N^*\cF)_\xi, (\pi^*E)_\xi)$
determine a well-defined section $\sigma_K$ of ${\mathcal
L}(\pi^*E)\otimes |T{\mathcal G}_N|^{1/2}$ over $G_{{\mathcal
F}_N}\cong N^*\cF \times_{T^*B} N^*\cF$. We say that the section $K$
is smooth, if the corresponding section $\sigma_K$ is smooth. This
defines an algebra isomorphism of $\Psi^{-\infty}(N^*\cF,\pi^*E)$
with $C^{\infty}(G_{{\mathcal F}_N}, {\mathcal L}(\pi^*E)\otimes
|T{\mathcal G}_N|^{1/2})$.

A Riemannian metric $g_M$ on $M$ is bundle-like if and only if there
exists a Riemannian metric $g_B$ on $B$ such that, for any $x\in M$,
the tangent map $dp$ induces an isometry from $(H_x, g_H)$ to
$(T_{p(x)}B, g_B)$, or, equivalently, $p : (M,g_M)\to (B,g_B)$ is a
Riemannian submersion. Then the transverse geodesic flow $f_t$ of
$g_M$ projects under the map $N^*\cF \to T^*B$ to the geodesic flow
$f^B_t$ of $g_B$ that implies commutativity of the following diagram
\[
  \begin{CD}
N^*\cF @>f_t>> N^*\cF\\ @VVV @VVV
\\ T^*B @>f^B_t>>
T^*B
  \end{CD}
\]
Commutativity of this diagram allows us to lift the flow $f_t$ to
the flow $F_t$ on the holonomy groupoid $G_{\cF_N}\cong
N^*\cF\times_{T^*B}N^*\cF$ in the following way:
\[
F_t(\nu_x,\nu_y)=(f_t(\nu_x),f_t(\nu_y)), \quad (\nu_x,\nu_y)\in
N^*\cF\times_{T^*B}N^*\cF.
\]

If $W$ is a Clifford bundle on $B$ equipped with a Clifford
connection $\nabla^W$ and $D_W$ is the associated Dirac operator
acting on $C^\infty(B,W)$, then its lift $p^*W$ to $M$ has a natural
structure of a transverse Clifford bundle, the pull-back of the
connection $\nabla^W$ is a leafwise flat, transverse Clifford
connection on $p^*W$, and the associated transverse Dirac operator
can be considered as a natural lift of $D_W$ to an operator in
$C^\infty(M,p^*W)$.

\begin{remark}
We refer the reader to \cite{matrix-egorov} for a discussion of a
particular case when the fibration $p: M\to B$ is a principal
$K$-bundle with a compact group $K$.
\end{remark}

\begin{remark}
Since any compact reduced orbifold is diffeomorphic to the quotient
of a compact manifold by an action of a compact Lie group with
finite isotropy groups, Theorems~\ref{Egorov1} and~\ref{Egorov2}
imply the Egorov theorem for Dirac-type operators on compact reduced
Riemannian orbifolds (see \cite{matrix-egorov} for more details).
\end{remark}

\subsection{Linear foliation on the two-torus}
Consider the linear foliation $\cF_\theta$ on the two-dimensional
torus $\TT^2=\RR^2/\ZZ^2$, whose leaves are the images of the lines
$\tilde{L}_{(u_0,v_0)}=\{(u_0+\tau, v_0+\theta \tau)\in \RR^2 :
\tau\in\RR\}, (u_0,v_0)\in\RR^2,$ under the projection $\RR^2\to
\TT^2$.

The conormal bundle $N^*\cF_\theta$ is described as
\[
N^*\cF_\theta=\{(u,v,p_u,p_v)\in T^*\TT^2\cong \TT^2\times \RR^2:
p_u+\theta p_v=0\},
\]
that is diffeomorphic to $\TT^2\times \RR$ with coordinates
$(u,v,p_v)$. The leaves of the lifted foliation $\cF_N$ are the
images of the lines $\tilde{L}_{(u_0,v_0,p^0_v)}=\{(u_0+\tau,
v_0+\theta \tau,p^0_v)\in \RR^2\times\RR : \tau\in\RR\},
(u_0,v_0,,p^0_v)\in\RR^2\times\RR,$ under the natural projection
$\RR^2\times\RR\to \TT^2\times\RR$.

If $\theta\in \QQ$, then all the leaves of $\cF_\theta$ are closed
and diffeomorphic to the circle $S^1$, and the foliation
$\cF_\theta$ is determined by the fibers of a fibration $\TT^2\to
S^1$. So this is a particular case of the situation discussed in
Section~\ref{s:Riem-sub}. The reduced phase space can be described
as $T^*S^1\cong S^1\times\RR$ with the coordinates $y=v-\theta u\in
S^1$ and $\eta=p_v\in \RR$, and the reduced symplectic structure is
given by the two-form $\omega=dy\wedge d\eta$.

Consider a Riemannian metric $g$ on $\TT^2$:
\[
g=a\, du^2+2b\,du\, dv +c\, dv^2,\quad a,b,c\in
C^\infty(\TT^2),\quad ac-b^2>0.
\]
Let $\omega_{\TT^2}$ be the Riemannian volume form on $\TT^2$:
\[
\omega_{\TT^2}=\sqrt{\det g_{\TT^2}}\,du\,dv. \quad \det g_{\TT^2}=
ac-b^2.
\]
The tangent space $T_{(u,v)}\cF_\theta\subset T_{(u,v)}\TT^2=\RR^2$
of $\cF_\theta$ is spanned by $(1,\theta)\in \RR^2$, and its
orthogonal complement $H_{(u,v)}$ is spanned by
$(-(b(u,v)+c(u,v)\theta), a(u,v)+b(u,v)\theta)\in \RR^2$. The
decomposition \eqref{e:g} holds with
\[
g_F=\frac{((a+b\theta)du+(b+c\theta)
dv)^2}{a+2b\theta+c\theta^2},\quad
g_H=\frac{ac-b^2}{a+2b\theta+c\theta^2}(dv-\theta du)^2.
\]
Therefore, the metric $g$ is bundle-like iff
\[
\left(\frac{\partial}{\partial u}+\theta\frac{\partial}{\partial
v}\right)a_H=0,
\]
where $a_H$ is the positive smooth function on $\TT^2$ defined by
\[
a^2_H=\frac{ac-b^2}{a+2b\theta+c\theta^2}.
\]
In particular, if $\theta$ is irrational, then $g$ is bundle-like
iff $a_H\equiv {\rm const}$.

There is an orthonormal base $\{e_F, e_H\}$ in $T \TT^2$ given by
\begin{align*}
e_F&=\frac{1}{\sqrt{a+2b\theta+c\theta^2}}\left(\frac{\partial}{\partial u}+\theta\frac{\partial}{\partial v}\right)\in F, \\
e_H&=\frac{1}{\sqrt{a+2b\theta+c\theta^2}}\frac{-(b+c\theta)\frac{\partial}{\partial
u}+ (a+b\theta)\frac{\partial}{\partial v}}{\sqrt{ac-b^2}}\in H.
\end{align*}
The dual orthonormal base $\{E_F,E_H\}$ in $T^*\TT^2$ is given by
\begin{align*}
E_F&=\frac{1}{\sqrt{a+2b\theta+c\theta^2}}\left((a+b\theta)du+
(b+c\theta)dv\right),
\\ E_H&=\frac{\sqrt{ac-b^2}}{\sqrt{a+2b\theta+c\theta^2}}(-\theta du+dv).
\end{align*}

The decomposition \eqref{e:d} holds with
\[
d_Ff=e_F(f)\cdot E_F,\quad  d_Hf=e_H(f)\cdot E_H, \quad \theta(f)=0,
\quad f\in C^\infty(\TT^2).
\]
The adjoint $d^*_H$ of the operator $d_H:C^\infty(\TT^2)\to
C^\infty(\TT^2, H^*)$ is given by
\[
d^*_H (g\cdot E_H) = e^*_H(g) = -e_H(g)+F\cdot g,
\]
where $F\in C^\infty(\TT^2)$ is given by
\begin{equation}\label{e:mean}
F=\frac{1}{\sqrt{ac-b^2}}\left(\frac{\partial }{\partial u}
\left(\frac{b+c\theta}{\sqrt{a+2b\theta+c\theta^2}}\right) -
\frac{\partial }{\partial v}
\left(\frac{a+b\theta}{\sqrt{a+2b\theta+c\theta^2}} \right)\right).
\end{equation}
The function $F$ is closely related with the mean curvature vector
$\tau$:
\[
\tau= Fe_H.
\]

For any $(u,v)\in \TT^2$, the complex Clifford algebra
$Cl(Q_{(u,v)})$ is the complex algebra generated by $1$ and $e_H$
with the relation $e_H^2=-1$. Thus, the vector bundle $Cl(Q)$ over
$\TT^2$ is trivial, and the sections $1\in C^\infty(\TT^2)$ and
$e_H\in C^\infty(\TT^2,Q)$ determine a natural trivialization of
this bundle. We also have $\nabla_{e_H}e_H =0$, and
\[
\nabla^{Cl(Q)}_{e_H} = e_H.
\]

Consider the transverse Clifford bundle $\cE$ defined in
Example~\ref{ex:sign}. Thus, we have $\cE=\Lambda H^*\otimes
\CC\cong \TT^2\times \CC^2$, It is equipped with the trivial
connection $\nabla^\cE$, and the action of $Cl(Q_{(u,v)})$ on
$\cE_{(u,v)}\cong \CC^2$ is given by
\[
c(e_H(u,v))=
\begin{pmatrix}
  0 & -1 \\
  1 & 0 \\
\end{pmatrix}
.
\]

The transverse Dirac operator $D_\cE$ associated with $\cE$ is a
matrix-valued first order differential operator acting on
$C^\infty(\TT^2,\TT^2\times \CC^2)$:
\begin{equation}\label{e:Dirac}
D_\cE=
\begin{pmatrix}
  0 & -e_H+\frac12 F \\
  e_H-\frac12 F & 0 \\
\end{pmatrix}
.
\end{equation}

The geodesic flow $f_t$ on $T^*\TT^2$ is the Hamiltonian flow of the
function
\[
h(u,v,p_u,p_v)=(G_{(u,v)}(p_u,p_v))^{1/2},
\]
where $G$ is the induced metric on $T^*\TT^2$ given by
\[
G_{(u,v)}(p_u,p_v)=\frac{1}{ac-b^2}(cp_u^2-2bp_up_v+ap_v^2), \quad
(p_u,p_v)\in T_{(u,v)}^*\TT^2.
\]

The infinitesimal generator $X$ of $f_t$ is given by
\begin{align*}
X=& (G_{(u,v)}(p_u,p_v))^{-1/2}
\frac{cp_u-bp_v}{ac-b^2}\frac{\partial}{\partial u}+
(G_{(u,v)}(p_u,p_v))^{-1/2}
\frac{-bp_u+ap_v}{ac-b^2}\frac{\partial}{\partial v}
\\ & -\frac{1}{2}(G_{(u,v)}(p_u,p_v))^{-1/2} \frac{\partial G_{(u,v)}(p_u,p_v)}{\partial u}
\frac{\partial}{\partial p_u}\\
& -\frac{1}{2}(G_{(u,v)}(p_u,p_v))^{-1/2} \frac{\partial
G_{(u,v)}(p_u,p_v)}{\partial v}\frac{\partial}{\partial p_v}.
\end{align*}

The restriction of $G$ to $N^*\cF_\theta$ is given by
\[
G_{(u,v)}(-\theta p_v,p_v)=\frac{1}{a_H^2}p^2_v, \quad p_v\in \RR.
\]
It follows that
\[
X(p_u+\theta p_v)=-\frac{1}{2}(G_{(u,v)}(p_u,p_v))^{-1/2}
\left(\frac{\partial }{\partial u}+\theta \frac{\partial }{\partial
v} \right) G_{(u,v)}(p_u,p_v),
\]
that immediately implies that, if the metric is bundle-like, then
the flow $f_t$ preserves $N^*\cF_\theta$. We also see that in this
case the restriction of $X$ to $N^*\cF_\theta$ is given by
\[
X=-a_H \frac{b+c\theta}{ac-b^2}{\rm sign}\,p_v
\frac{\partial}{\partial u}+a_H \frac{a+b\theta}{ac-b^2}{\rm
sign}\,p_v \frac{\partial}{\partial v}+ a^{-2}_H \frac{\partial
a_H}{\partial v} |p_v|\frac{\partial}{\partial p_v}.
\]
We have
\[
X(v-\theta u)=-\frac{1}{a_H}{\rm sign}\,p_v,
\]
that is constant along the leaves of $\cF_N$, because $a_H$ depends
only on $v-\theta u$. It follows that the flow $f_t$ takes each leaf
of $\cF_N$ to a leaf.

So if $\theta\in \QQ$, the reduced dynamic on the reduced phase
space $T^*S^1$ with coordinates $y=v-\theta u$ and $\eta=p_v$ is
determined by the reduced Hamiltonian $h(y,\eta)
=\frac{1}{a_H(y)}|\eta|$:
\[
\frac{dy}{dt}=-\frac{1}{a_H(y)}{\rm sign}\,\eta, \quad
\frac{d\eta}{dt}= a^{-2}_H \frac{\partial a_H}{\partial y} |\eta|.
\]

If $\theta\not\in \QQ$, the foliation $\cF_\theta$ is given by the
orbits of a free action of $\RR$ on $\TT^2$, and its holonomy
groupoid is described as follows: $G=\TT^2\times \RR$,
$G^{(0)}=\TT^2$, $s(u,v,\tau)=(u-t,v-\theta \tau)$,
$r(u,v,\tau)=(u,v)$, $(u,v)\in \TT^2, \tau\in\RR,$ and the product
of $(u_1,v_1,\tau_1)$ and $(u_2,v_2,\tau_2)$ is defined if
$u_2=u_1-\tau_1, v_2=v_1-\theta \tau_1$ and equals
\[
(u_1,v_1,\tau_1)(u_2,v_2,\tau_2)=(u_1,v_1,\tau_1+\tau_2).
\]

\begin{remark}
If $\theta\in\QQ$, then the linear foliation on $\TT^2$ is given by
the orbits of a free group action of $S^1$ on $\TT^2$, and its
holonomy groupoid coincides with the crossed product groupoid
$G=\TT^2\rtimes S^1$.
\end{remark}

The holonomy groupoid $G_{\cF_N}$ is described as follows:
$G_{\cF_N}=\TT^2\times \RR^2$, $G_{\cF_N}^{(0)}=\TT^2\times \RR$,
$s_N(u,v,p_v,\tau)=(u-t,v-\theta \tau,p_v)$,
$r_N(u,v,p_v,\tau)=(u,v,p_v)$, $(u,v,p_v)\in \TT^2\times \RR,
\tau\in\RR,$ and the product of $(u_1,v_1,p_{v,1},\tau_1)$ and
$(u_2,v_2,p_{v,2},\tau_2)$ is defined if $u_2=u_1-\tau_1,
v_2=v_1-\theta \tau_1, p_{v,1}=p_{v,2}(=p_v),$ and equals
\[
(u_1,v_1,p_v,\tau_1)(u_2,v_2,p_v,\tau_2)=(u_1,v_1,p_v,\tau_1+\tau_2).
\]

There is a natural trivialization of the line bundle $T{\mathcal
F}_N$ denoted by $d\tau$. Let $\omega_F\in C^\infty
(\TT^2,T{\mathcal F}_N)$ be the leafwise Riemannian volume form:
\[
\omega_F(u,v,p_v)=\sqrt{\det g_F(u,v)}d\tau, \quad \det
g_F=a+2b\theta+c\theta^2.
\]
An arbitrary section $\sigma \in C^{\infty}(G_{{\mathcal F}_N},
|T{\mathcal G}_N|^{1/2})$ can be written as
\begin{multline}\label{e:defsigma}
\sigma(u,v,p_v,\tau) = k(u,v,p_v,\tau)|\omega_F(u,v,p_v)|^{1/2}
|\omega_F(u-\tau,v-\theta \tau,p_v)|^{1/2}, \\ (u,v)\in \TT^2,\quad
p_v\in \RR,\quad \tau\in\RR,
\end{multline}
It is properly supported if and only if $k$ has compact support as a
function of $\tau$ for any $(u,v,p_v)$.

For any $\sigma_j\in C^{\infty}_{prop}(G_{{\mathcal F}_N},
|T{\mathcal G}_N|^{1/2}),$ of the form \eqref{e:defsigma} with some
$k_j, j=1,2$, their product $\sigma=\sigma_1\ast \sigma_2$ is
written in the form \eqref{e:defsigma} with
\begin{multline*}
(k_1\ast k_2)(u,v,p_v,\tau)\\
=\int_{-\infty}^\infty k_1(u,v,p_v,\tau_1)k_2(u-\tau_1, v-\theta
\tau_1, p_v, \tau-\tau_1)\, \sqrt{\det g_F(u-\tau_1,v-\theta
\tau_1)}\,d\tau_1, \\
(u,v)\in \TT^2,\quad p_v \in \RR, \quad \tau\in\RR.
\end{multline*}

The infinitesimal generator $\cH$ of the induced flow $F_t$ on
$G_{{\mathcal F}_N}$ is given by
\begin{equation*}
\begin{aligned}
\cH=& -a_H \frac{b(u,v)+c(u,v)\theta}{a(u,v)c(u,v)-b(u,v)^2}{\rm sign}\,p_v \frac{\partial}{\partial u} \\
& + a_H \frac{a(u,v)+b(u,v)\theta}{a(u,v)c(u,v)-b(u,v)^2}{\rm
sign}\,p_v\frac{\partial}{\partial v}\\ & +
a_H\Big(\frac{b(u-\tau,v-\theta\tau)+c(u-\tau,v-\theta\tau)\theta}{a(u-\tau,v-\theta\tau)
c(u-\tau,v-\theta\tau)-b(u-\tau,v-\theta\tau)^2} \\ & -
\frac{b(u,v)+c(u,v)\theta}{a(u,v)c(u,v)-b(u,v)^2}\Big){\rm
sign}\,p_v\frac{\partial}{\partial \tau}.
\end{aligned}
\end{equation*}

We have
\[
\cL_{X}|\omega_F|^{1/2}(u,v,p_v)=\frac12
F(u,v)|\omega_F|^{1/2}(u,v,p_v),
\]
where the mean curvature $F$ is given by \eqref{e:mean}. Thus, for
any $\sigma\in C^{\infty}(G_{{\mathcal F}_N}, |T{\mathcal
G}_N|^{1/2})$ of the form \eqref{e:defsigma}, we have
\begin{multline}\label{e:defH}
\cL_{\cH}\sigma(u,v,p_v,\tau) = \left(\cH +\frac12 F(u,v)+\frac12
F(u-\tau,v-\theta \tau)\right)k(u,v,p_v,\tau)\times \\ \times
|\omega_F(u,v,p_v)|^{1/2} |\omega_F(u-\tau,v-\theta
\tau,p_v)|^{1/2}, \quad (u,v)\in \TT^2,\quad p_v\in \RR,\quad
\tau\in\RR,
\end{multline}

The bundle $\pi^*\cE$ is a trivial two-dimensional complex bundle
over $N^*\cF_\theta \cong \TT^2\times \RR$. Therefore, an element
$\sigma\in C^{\infty}(G_{{\mathcal F}_N},{\mathcal
L}(\pi^*\cE)\otimes |T{\mathcal G}_N|^{1/2})$ can be written as
\begin{multline*}
\sigma(u,v,p_v,\tau) = k(u,v,p_v,\tau)|\omega_F(u,v,p_v)|^{1/2}
|\omega_F(u-\tau,v-\theta \tau,p_v)|^{1/2}, \\ (u,v,p_v,\tau)\in
G_{{\mathcal F}_N}\cong \TT^2 \times \RR \times \RR,
\end{multline*}
where $k$ is a smooth function with values in the space $M_2(\CC)$
of complex $2\times 2$ matrices. Since the Clifford connection
$\nabla^\cE$ is trivial, the action of $\cH$ on such a $\sigma$ is
defined by the formula \eqref{e:defH} with
\[
(\cH k)_{\alpha\beta}=\cH (k_{\alpha\beta}), \quad \alpha,\beta=1,2.
\]

The class $\Psi^{d,-\infty}(\TT^2,{\mathcal F}_\theta,\TT^2\times
\CC^2)$, $d\in \RR$, can be described as follows. Suppose that a
polyhomogeneous symbol $k\in S^d(\RR^4\times \RR, M_2(\CC))$
satisfies the conditions:
\begin{enumerate}
  \item for any $(m,n)\in\ZZ^2$,
\begin{equation}\label{e:K1}
k(u+m,v+n,u'+m,v'+n,\eta)=k(u,v,u',v',\eta),\quad
(u,v,u',v',\eta)\in \RR^4\times \RR;
\end{equation}
  \item there exists $R>0$ such that
\begin{equation}\label{e:K2}
k(u,v,u',v',\eta)=0,\quad (u,v,u',v',\eta)\in \RR^4\times \RR, \quad
(u-u')^2+(v-v')^2>R^2.
\end{equation}
\end{enumerate}

Define an operator $K:C^\infty(\RR^2,\RR^2\times\CC^2)\to
C^\infty(\RR^2,\RR^2\times\CC^2)$ by the formula
\begin{multline}\label{e:defK}
Kf(u,v)\\ =\int e^{i(v-v'-\theta(u-u'))\eta}\, k(u,v,u',v',\eta)
\,f(u',v')\,\sqrt{\det g_F(u^\prime,v^\prime)}
\,du^\prime\,dv^\prime\,d\eta.
\end{multline}
By \eqref{e:K1} and \eqref{e:K2}, the operator $K$ takes any
$\ZZ^2$-invariant function from $C^\infty(\RR^2,\RR^2\times\CC^2)$
to a $\ZZ^2$-invariant function from
$C^\infty(\RR^2,\RR^2\times\CC^2)$, and determines an operator,
acting on $C^\infty(\TT^2,\TT^2\times\CC^2)$.

It is clear that any operator from $\Psi^{d,-\infty}(\TT^2,{\mathcal
F}_\theta,\TT^2\times\CC^2)$ can be written in this form. On the
other hand, one can show that any operator $K$ defined by
\eqref{e:defK} belongs to $\Psi^{d,-\infty}(\TT^2,{\mathcal
F}_\theta,\TT^2\times\CC^2)$, and its principal symbol is given by
\begin{multline*}
\sigma(K)(u,v,p_v,\tau) = \sum_{(m,n)\in \ZZ^2}
k_d(u+m,v+n,u-\tau,v-\theta \tau,p_v)\times
\\ \times |\omega_F(u,v,p_v)|^{1/2} |\omega_F(u-\tau,v-\theta
\tau,p_v)|^{1/2}, \quad (u,v,p_v,\tau)\in G_{{\mathcal F}_N}\cong
\TT^2 \times \RR \times \RR,
\end{multline*}
where $k_d$ is the (degree $d$) homogeneous component of $k$.

The orbits of the vector field $e_H$ determine a foliation,
$\cF^\pitchfork$, transverse to $\cF_\theta$. It is easy to see
that, for any natural $m$, the space $\cH^m$ coincides with the
space of all $f\in L^2(\TT^2)$ such that $e_H^jf\in L^2(\TT^2)$ for
$j=1,\ldots,m$. Thus, for any $s\in \RR$, the space $\cH^s$
coincides with the anisotropic Sobolev space $H^{0,s}(\TT^2,
\cF^\pitchfork)$ \cite{tang}.

Now Theorems~\ref{Egorov1} and~\ref{Egorov2} read as follows.

\begin{theorem}
\label{Egorov-torus} Let $D_\cE$ be the transverse Dirac operator
defined by~\eqref{e:Dirac}.
\medskip\par
(1) For any $K\in \Psi^{m,-\infty}(\TT^2,{\mathcal
F}_\theta,\TT^2\times\CC^2)$, there exists an operator
$K(t)\in\Psi^{m,-\infty}(\TT^2,{\mathcal
F}_\theta,\TT^2\times\CC^2)$ such that $\Phi_t(K)-K(t), t\in \RR,$
is a smooth family of operators of class
$\cL^1(\cH^{-\infty},\cH^{\infty})$.
\medskip\par
(2) If $\sigma \in S^m(G_{{\mathcal F}_N},{\mathcal
L}(\pi^*\cE)\otimes |T{\mathcal G}_N|^{1/2})$ is the principal
symbol of $K$:
\[
\sigma(u,v,p_v,\tau) = k(u,v,p_v,\tau)|\omega_F(u,v,p_v)|^{1/2}
|\omega_F(u-\tau,v-\theta \tau,p_v)|^{1/2},
\]
then the principal symbol $\sigma_t\in S^m(G_{{\mathcal
F}_N},{\mathcal L}(\pi^*\cE)\otimes |T{\mathcal G}_N|^{1/2})$ of
$K(t)$:
\[
\sigma_t(u,v,p_v,\tau) = k(t,u,v,p_v,\tau)|\omega_F(u,v,p_v)|^{1/2}
|\omega_F(u-\tau,v-\theta \tau,p_v)|^{1/2},
\]
is the solution of the equation
\[
\frac{\partial k}{\partial t}(t,u,v,p_v,\tau) = \left(\cH+\frac12
F(u,v)+\frac12 F(u-\tau,v-\theta \tau)\right)k(t,u,v,p_v,\tau),
\]
satisfying the condition $k(0,u,v,p_v,\tau)=k(u,v,p_v,\tau)$.
\end{theorem}

\bibliographystyle{amsplain}


\providecommand{\bysame}{\leavevmode\hbox
to3em{\hrulefill}\thinspace}
\providecommand{\MR}{\relax\ifhmode\unskip\space\fi MR }
\providecommand{\MRhref}[2]{%
  \href{http://www.ams.org/mathscinet-getitem?mr=#1}{#2}
} \providecommand{\href}[2]{#2}


\begin{thebibliography}{10}

\bibitem{Block-Ge}
Jonathan Block and Ezra Getzler, \emph{Quantization of foliations},
Proceedings
  of the {XX}th {I}nternational {C}onference on {D}ifferential {G}eometric
  {M}ethods in {T}heoretical {P}hysics, {V}ol.\ 1, 2 ({N}ew {Y}ork, 1991),
  World Sci. Publ., River Edge, NJ, 1992, pp.~471--487. \MR{1225136
  (95d:58149)}

\bibitem{Co79}
Alain Connes, \emph{Sur la th\'eorie non commutative de
l'int\'egration},
  Alg\`ebres d'op\'erateurs ({S}\'em., {L}es {P}lans-sur-{B}ex, 1978), Lecture
  Notes in Math., vol. 725, Springer, Berlin, 1979, pp.~19--143. \MR{548112
  (81g:46090)}

\bibitem{Co}
\bysame, \emph{Noncommutative geometry}, Academic Press Inc., San
Diego, CA,
  1994. \MR{1303779 (95j:46063)}

\bibitem{G-K91a}
James~F. Glazebrook and Franz~W. Kamber, \emph{On spectral flow of
transversal
  {D}irac operators and a theorem of {V}afa-{W}itten}, Ann. Global Anal. Geom.
  \textbf{9} (1991), no.~1, 27--35. \MR{1116629 (92i:58180)}

\bibitem{G-K91}
\bysame, \emph{Transversal {D}irac families in {R}iemannian
foliations}, Comm.
  Math. Phys. \textbf{140} (1991), no.~2, 217--240. \MR{1124268 (92j:58103)}

\bibitem{GS79}
Victor Guillemin and Shlomo Sternberg, \emph{Some problems in
integral geometry
  and some related problems in microlocal analysis}, Amer. J. Math.
  \textbf{101} (1979), no.~4, 915--955. \MR{536046 (82b:58087)}

\bibitem{Gu-S82}
\bysame, \emph{Homogeneous quantization and
  multiplicities of group representations}, J. Funct. Anal. \textbf{47} (1982),
  no.~3, 344--380. \MR{665022 (84d:58034)}

\bibitem{Gu-Uribe}
Victor Guillemin and Alejandro Uribe, \emph{Circular symmetry and
the trace formula},  Invent. Math. \textbf{96} (1989), no.~2,
385--423. \MR{989702 (90e:58159)}

\bibitem{Gu-Uribe90}
\bysame, \emph{Reduction and the trace formula}, J. Differential
Geom.  \textbf{32} (1990), no.~2, 315--347. \MR{1072909 (92h:58191)}

\bibitem{tang}
Yuri~A. Kordyukov, \emph{Functional calculus for tangentially
elliptic
  operators on foliated manifolds}, Analysis and geometry in foliated manifolds
  ({S}antiago de {C}ompostela, 1994), World Sci. Publ., River Edge, NJ, 1995,
  pp.~113--136. \MR{1414199 (97h:58158)}

\bibitem{noncom}
\bysame, \emph{Noncommutative spectral geometry of {R}iemannian
foliations},
  Manuscripta Math. \textbf{94} (1997), no.~1, 45--73. \MR{1468934
  (98j:58116)}

\bibitem{egorgeo}
\bysame, \emph{Egorov's theorem for transversally elliptic operators
on
  foliated manifolds and noncommutative geodesic flow}, Math. Phys. Anal. Geom.
  \textbf{8} (2005), no.~2, 97--119. \MR{2160331 (2006e:58038)}

\bibitem{matrix-egorov}
\bysame, \emph{The {E}gorov theorem for transverse {D}irac-type
operators on
  foliated manifolds}, J. Geom. Phys. \textbf{57} (2007), no.~11, 2345--2364.
  \MR{2360245 (2008k:58061)}

\bibitem{survey}
\bysame, \emph{Noncommutative geometry of foliations}, J. K-Theory
\textbf{2}
  (2008), no.~2, Special issue in memory of Yurii Petrovich Solovyev. Part 1,
  219--327. \MR{2456103}

\bibitem{vanishing}
\bysame, \emph{Vanishing theorem for transverse {D}irac operators on
  {R}iemannian foliations}, Ann. Global Anal. Geom. \textbf{34} (2008), no.~2,
  195--211. \MR{2425530}

\bibitem{LM87}
Paulette Libermann and Charles-Michel Marle, \emph{Symplectic
geometry and
  analytical mechanics}, Mathematics and its Applications, vol.~35, D. Reidel
  Publishing Co., Dordrecht, 1987, Translated from the French by Bertram Eugene
  Schwarzbach. \MR{882548 (88c:58016)}

\bibitem{Li75}
Andr{\'e} Lichnerowicz, \emph{Vari\'et\'e symplectique et dynamique
associ\'ee
  \`a une sous-vari\'et\'e}, C. R. Acad. Sci. Paris S\'er. A-B \textbf{280}
  (1975), A523--A527. \MR{0405503 (53 \#9296)}

\bibitem{Li77}
\bysame, \emph{Les vari\'et\'es de {P}oisson et leurs alg\`ebres de
{L}ie
  associ\'ees}, J. Differential Geometry \textbf{12} (1977), no.~2, 253--300.
  \MR{0501133 (58 \#18565)}

\bibitem{MW}
Jerrold Marsden and Alan Weinstein, \emph{Reduction of symplectic
manifolds
  with symmetry}, Rep. Mathematical Phys. \textbf{5} (1974), no.~1, 121--130.
  \MR{0402819 (53 \#6633)}

\bibitem{Molino}
Pierre Molino, \emph{Riemannian foliations}, Progress in
Mathematics, vol.~73,
  Birkh\"auser Boston Inc., Boston, MA, 1988, Translated from the French by
  Grant Cairns, With appendices by Cairns, Y. Carri{\`e}re, {\'E}. Ghys, E.
  Salem and V. Sergiescu. \MR{932463 (89b:53054)}

\bibitem{Re}
Bruce~L. Reinhart, \emph{Differential geometry of foliations},
Ergebnisse der
  Mathematik und ihrer Grenzgebiete,
  vol.~99, Springer-Verlag, Berlin, 1983, The fundamental integrability
  problem. \MR{705126 (85i:53038)}

\bibitem{Tate99}
Tatsuya Tate, \emph{Quantum ergodicity at a finite energy level}, J.
Math. Soc.  Japan \textbf{51} (1999), no.~4, 867--885. \MR{1705252
(2001f:81051)}

\bibitem{Xu}
Ping Xu, \emph{Noncommutative {P}oisson algebras}, Amer. J. Math.
\textbf{116}
  (1994), no.~1, 101--125. \MR{1262428 (94m:58246)}

\bibitem{Zelditch92}
Steven Zelditch, \emph{On a ``quantum chaos'' theorem of {R}.
{S}chrader and  {M}. {T}aylor}, J. Funct. Anal. \textbf{109} (1992),
no.~1, 1--21.  \MR{1183602 (94d:58155)}

\end{thebibliography}
\end{document}